# SOME APPLICATIONS OF THE KOHN–ROSSI EXTENSION THEOREMS

GEORGE MARINESCU

ABSTRACT. We prove extension results for meromorphic functions by combining the Kohn-Rossi extension theorems with Andreotti's theory on the algebraic and analytic dependence of meromorphic functions on pseudoconcave manifolds. Versions of Kohn-Rossi theorems for pseudoconvex domains are included.

## 1. INTRODUCTION

A naive idea of extending a meromorphic function $\mu$ beyond its domain of definition is to write $\mu$ as a quotient of two holomorphic functions and then extend these functions. At least for holomorphic functions defined near the boundary of complex manifolds this idea works as in the proof of the Kohn-Rossi extension theorem (see Theorem 2.2 and the proof of subsequent Theorem 3.3). But if $\mu$ is defined in a neighbourhood $U$ of the complement $Y \setminus M$ of a smooth domain $M$ whose boundary Levi form has at least one positive eigenvalue in a compact complex manifold $Y$, we cannot apply this method. Indeed, $U$ is pseudoconcave and any holomorphic function on $U$ is locally constant. However, it is rather simple to see that a meromorphic function on any complex manifold can be written as a quotient of two holomorphic sections of an appropriate holomorphic line bundle. We can thus hope to extend the line bundle and its sections. The situation in the above example becomes particulary favourable in this respect if we assume $Y$ to be Moishezon. Then by Andreotti's theory [1] it can be shown that $\mu$ can be written as a quotient of two holomorphic sections in some tensor power $E^m$ of a line bundle $E \to Y$ with maximum Kodaira-Iitaka dimension. Therefore we need to extend but the sections of $E^m$ and this is done using the theorems of Kohn-Rossi. The situation being global we have to impose some supplementary hypotheses, namely one of the following (see Theorem 2.4):

(H1) $M$ is $q$–complete ($q < n$).
(H2) There exists one non-constant holomorphic function in the neighbourhood of $bM$.

1991 *Mathematics Subject Classification.* Primary 32D15, 32A20, Secondary 32F10.
The author was supported by a NATO/Royal Society Fellowship at the University of Edinburgh. He thanks the Department of Mathematics & Statistics for its hospitality.

(H3) The holomorphic line bundle $E$ of maximal Kodaira-Iitaka dimension is trivial in the neighbourhood of $M$.

Note that using the Thullen-type theorem of Siu [20] (see also Ivashkovich [11] §4 and Shiffman [21]) about the extension of a meromorphic map from the Hartogs figure to a polydisc one can prove that if $\mu$ is a meromorphic map defined in $Y \setminus \overline{M}$ with values in a Kähler manifold then $\mu$ has a meromorphic continuation into $bM$. Our result can be viewed as global counterpart of the previous theorem. As for the best global result let us mention the following theorem of Ivashkovich [11]: every meromorphic map from a domain $D$ in a Stein manifold into a compact Kähler manifold $K$ extends to a meromorphic map from the envelope of holomorphy $\widehat{D}$ of $D$ into $K$.

The Kohn-Rossi theorems are consequences of the global regularity of the $\bar{\partial}$–Neumann problem for domains $M$ as above. Since the global regularity holds also for some class of pseudoconvex domains [14], [17] we prove in the last section versions of those theorems for this class of domains.

*Aknowledgements.* I wish to acknowledge my debt to Professor Guennadi Henkin, who is at the origin of the ideas used in this article. He deserves a large share of the credit for any originality the present paper may possess.

## 2. Meromorphic continuation

Let us give a brief account of the extension theory developed by Kohn and Rossi. Let $M$ be a smooth domain in a complex manifold. A smooth function on $bM$ is said to be Cauchy-Riemann ($CR$) if its differential vanishes on the sub-bundle $\mathbb{C}T(bM) \cap T^{0,1}M$. A smooth section over $bM$ of a holomorphic vector bundle is said to be $CR$ if it is represented in local holomorphic frames by vectors of $CR$ functions. In both cases we write $\bar{\partial}_b f = 0$. Of course, the restriction of a holomorphic section to $bM$ is $CR$. The crucial lemma depends on the natural duality between the free boundary conditions imposed by $\bar{\partial}^*$ and $\vartheta^*$ ($\vartheta$ is the formal adjoint and $\bar{\partial}^*$ is the hilbertian adjoint of $\bar{\partial}$) and the solution of the $\bar{\partial}$–Neumann problem. To simplify matters let us take $\varphi \in \mathcal{E}^{0,0}(V, F)$ merely a holomorphic section of the holomorphic vector bundle $F$ over a neighbourhood of $bM$. We define the Hodge operator $\# : \mathcal{E}^{p,q}(F) \longrightarrow \mathcal{E}^{n-p,n-q}(F^*)$ by $\alpha \wedge \#\beta = <\alpha, \beta> d\,vol$, where the exterior product $\wedge$ is combined with the canonical pairing $F \times F^* \to \mathbb{C}$. Recall that $\vartheta = -\#\bar{\partial}\#$ and $\#\bar{\partial} = (-1)^{r+1}\vartheta\#$ on $r$–forms. Recall also that if a form with $\mathcal{C}^1$ coefficients on $\overline{M}$ vanishes on $bM$ then clearly belongs to $\mathrm{Dom}(\bar{\partial}^*)$ (integration by parts) and $\bar{\partial}^* = \vartheta$ on such forms. Consider $\varphi'$ a smooth extension of $\varphi$ to $M$ and assume that the $\bar{\partial}$–Neumann problem is solvable in bi-degree $(n, n-1)$ for $F^*$ (which is the case if the Levi form of $bM$ has one positive eigenvalue).

Denote by $N : \mathcal{E}^{n,n-1}(\overline{M}, F^*) \longrightarrow \mathcal{E}^{n,n-1}(\overline{M}, F^*)$ the Neumann operator. Define
$$\varphi_0 = \varphi' + \#\bar{\partial}N\#\bar{\partial}\varphi' \in \mathcal{E}^{0,0}(\overline{M}, F).$$

Consider now $\alpha = \#\bar\partial\varphi' \in \mathcal{E}^{n,n-1}(\overline{M}, F^*)$ which has compact support in $M$ (since $\bar\partial\varphi'$ has). Clearly $\alpha \in \text{Dom}(\bar\partial^*)$ and $\alpha = -\vartheta\#\varphi'$. (If $\varphi$ is only $CR$ then of course $\bar\partial\varphi'$ hasn't compact support, but we can construct an extension $\varphi'$ such that $\bar\partial\varphi'|_{bM} = 0$; alternatively, Lemma 2.9 of [13] shows that in fact $\bar\partial_b\varphi = 0$ is equivalent to $\alpha = \vartheta\#\varphi' \in \text{Dom}(\bar\partial^*)$ for any extension $\varphi'$.) Thus,

$$\bar\partial^* N\alpha = N\bar\partial^*\alpha = N\bar\partial^*\vartheta\#\varphi' = 0,$$

so that the Hodge decomposition given by the solution of the $\bar\partial$–Neumann $\alpha = \bar\partial\bar\partial^* N\alpha + \bar\partial^*\bar\partial N\alpha + H\alpha$ reduces to $\alpha = \bar\partial^*\bar\partial N\alpha + H\alpha$ where $H$ is the orthogonal projector on the harmonic space $\mathcal{H}^{n,n-1}(M, F^*)$, the space of smooth forms $\chi \in \text{Dom}(\bar\partial^*)$ such that $\bar\partial\chi = 0$ and $\bar\partial^*\chi = 0$. Assume that

(2.1) $$\alpha \perp \mathcal{H}^{n,n-1}(M, F^*).$$

Then if we denote $\psi = \bar\partial N\alpha$ we have $\alpha = \bar\partial^*\psi$ and $\#\varphi_0 = \#\varphi' + \psi$. Therefore, $\vartheta\#\varphi_0 = \vartheta\#\varphi' + \vartheta\psi = -\alpha + \vartheta\psi = 0$, which implies $\#\bar\partial\varphi_0 = 0$, $\bar\partial\varphi_0 = 0$. On the other hand we have $\psi \in \text{Dom}(\bar\partial^*) \cap \mathcal{E}^{n,n}(\overline{M}, F^*)$. In this extreme case the boundary conditions imposed by the preceding relation are Dirichlet conditions, $\psi|_{bM} = 0$ so that $\varphi_0 = \varphi' = \varphi$ on $bM$. We easily infer that $\varphi_0 = \varphi$ in a neighbourhood of $bM$ since real hypersurfaces are uniqueness sets for holomorphic functions and sections in holomorphic vector bundles.

Conversely, if $\varphi$ admits a holomorphic extension $\varphi'$ then $\alpha$ and $H\alpha$ vanish, so that condition (2.1) is necessary and sufficient for the holomorphic extension. Integrating by parts

$$(\chi, \alpha) = (\chi, \vartheta\#\varphi') = (\bar\partial\chi, \#\varphi') + \int_{bM} \chi \wedge \varphi$$

for $\chi \in \mathcal{H}^{n,n-1}(M, F^*)$ we obtain:

**Lemma 2.1.** *Let $M$ be a relatively compact domain with smooth boundary in a complex manifold $Y$ such that the Levi form of the boundary $bM$ has at least one positive eigenvalue everywhere. Let $F$ be a holomorphic line bundle over $Y$. Then a given smooth CR section $\sigma$ of $F$ over $bM$ has a smooth extension over $\overline{M}$ which is holomorphic in $M$ if and only if*

(2.2) $$\int_{bM} \chi \wedge \sigma = 0 \quad \text{for any} \quad \chi \in \mathcal{H}^{n,n-1}(M, F^*)$$

On grounds of finite dimensionality of the harmonic space [13] the Lemma entails (see also the proof of Theorem 3.3):

**Theorem 2.2** (Kohn-Rossi). *Let $M$ as in Lemma 2.1 and with connected boundary. Then every CR function on $bM$ has a holomorphic extension to all of $M$.*

Note that if there exist a non-constant $CR$ functions connectedness of $bM$ is a superfluous hypothesis.

**Theorem 2.3** (Kohn-Rossi). *Let $Y$ be a complex manifold possesing at least one non-constant global holomorphic function. Let $M$ as in Lemma 2.1.*

*Then for every holomorphic line bundle $E$ on $Y$ there exists a holomorphic line bundle $L$ suct that every smooth* CR *section of $E \otimes L$ over the boundary $bM$ extends by a smooth section of $E \otimes L$ over $\overline{M}$, holomorphic in $M$.*

Of course these three statements remain true for holomorphic sections defined near $bM$ thanks to the uniqueness theorem. Relying on these results we shall prove the following:

**Theorem 2.4.** *Let $Y$ be a Moishezon manifold $Y$ of dimension $n \geqslant 2$. Assume that $M$ is a smooth domain in $Y$ subject to one of the following hypothesis:*

- **(H1):** $M$ *is $q$–complete ($q < n$).*
- **(H2):** *The Levi form of $bM$ has one positive eigenvalue everywhere and there exists one non-constant holomorphic function in the neighbourhood of $bM$.*
- **(H3):** *The Levi form of $bM$ has one positive eigenvalue everywhere and there exists a holomorphic line bundle $E$ on $Y$ of maximal Kodaira-Iitaka dimension which is trivial in the neighbourhood of $M$.*

*Then any meromorphic function defined in a pseudoconcave neighbourhood $U$ of $bM$ (e.g., a connected neighbourhood of $Y \setminus M$), there exists a meromorphic function on $M \cup U$ which agrees with it on $U$.*

Observe that in view of the theorem of Siu and Ivashkovich quoted in the Introduction we may assume that $\mu$ is defined only on $U \cap (Y \setminus \overline{M})$.

Recall that a compact connected complex manifold is said to be *Moishezon* if the transcendence degree of its meromorphic function field equals its complex dimension. The *Kodaira-Iitaka dimension* $\kappa(E)$ of a holomorphic line bundle on a compact complex manifold is the supremum of the generic rank of the canonical rational maps defined by the nonzero sections in $H^0(E^m)$ for $m \geqslant 1$ (if any) and $\kappa(E) = -\infty$ otherwise. By a well-known result of Siegel we have $0 \leqslant \kappa(E) \leqslant \deg \operatorname{tr} \mathcal{K}(Y) \leqslant n = \dim Y$ where $\mathcal{K}(Y)$ is the field of meromorphic functions on $Y$. Therefore, if $Y$ possess a line bundle with maximum Kodaira-Iitaka dimension then it is automatically Moishezon. Conversely, if $Y$ is Moishezon then there exists a nonsingular projective modification (composition of blow-ups of nonsingular centres) $\pi : \widetilde{Y} \longrightarrow Y$; if we consider an ample divisor $\widetilde{H}$ on $\widetilde{Y}$ then the line bundle associated to $H = \pi_* \widetilde{H}$ has maximum Kodaira-Iitaka dimension.

Let us introduce now the convexity notions used. A real smooth function $\varrho$ defined on a complex manifold $X$ of dimension $n$ is said to be $q$–convex if its complex hessian $i \partial \bar{\partial} \varrho$ has at least $n - q + 1$ positive eigenvalues at every point of $X$; $\varrho$ is called exhaustion function if for every $c \in \mathbb{R}$ the set $\{\varrho < c\}$ is relatively compact. A complex manifold is said to be $q$–*complete* if it is endowed with a $q$–convex exhaustion function. Let us consider now a relatively compact smooth domain $M$ whose boundary has a definition function $\varrho$; the Levi form of $bM$ is the complex hessian of $\varrho$ restricted to the holomorphic tangent space $T^{1,0}M \cap \mathbb{C}T(bM)$.

Following Andreotti [1], a connected complex $U$ manifold is called *pseudoconcave* if one can find a non-empty open subset $N \Subset U$ with psudoconcave boundary $bN$, i.e., for any point $z \in bN$, given a holomorphic function defined in the neighbourhood of $z$, one can find a neighbourhood $V$ of $z$ such that maximum of $|f|$ on $V \cap \overline{N}$ is attained in a point of $V \cap N$. It is well-known [1] that if $Y$ is compact and $M \subset Y$ is a smooth domain whose boundary Levi form has one positive eigenvalue everywhere then any neighbourhood $U$ of $Y \setminus M$ is pseudoconcave. Also, any compact manifold is pseudoconcave: take $N$ to be the complement of a twisted ball (a stongly pseudoconvex set – the Levi form of the boundary is positive definite).

Let us consider a connected pseudoconcave manifold $U$ of dimension $n$ and a holomorphic line bundle $E$ over $U$. We define the graded ring

$$\mathcal{A}(E) = \bigoplus_{m \geqslant 0} H^0(X, E^m).$$

Since $U$ is connected $\mathcal{A}(E)$ is an integral domain and we can consider the quotient field $\mathcal{Q}(E)$; $U$ being a manifold (normal space), $\mathcal{Q}(E)$ is algebraically closed in the meromorphic function field $\mathcal{K}(U)$.

The vector space $H^0(X, E^m)$ is finite dimensional and if $s_0, \ldots, s_k$ is a basis we define the canonical rational maps

$$\Phi_m : X -\!-\!\to \mathbb{P}^{k-1}, \quad x \rightsquigarrow [s_0(x), \ldots, s_k(x)].$$

We say that $\mathcal{A}(E)$ gives local coordinates at a point $x$ if there exists $m \geqslant 1$ such that $\Phi_m$ is holomorphic at $x$ and $\operatorname{rank}_x \Phi_m = n$. If $U$ is compact then $\mathcal{A}(E)$ gives local coordinates at some point of $U$ if and only if $\kappa(E) = n$; in that case $\mathcal{A}(E)$ gives local coordinates on an open dense set.

*Proof of Theorem 2.4.* We use the following basic result of Andreotti [1] which asserts that on a pseudoconcave manifold $U$ the algebraic dependence of meromorphic functions is the same as the analytic dependence. A corollary of this is that the fields $\mathcal{K}(U)$ and $\mathcal{Q}(F)$, for some holomorphic line bundle on $U$, are isomorphic to simple algebraic extensions of fields of rational functions with $d$ and $e$ variables, respectively, where $e \leqslant d \leqslant n$.

Let us consider now the line bundle $E$ given by the hypothesis. Since $\kappa(E) = n$, $\mathcal{A}(E)$ gives local coordinates on an open dense subset of $Y$, and a fortiori of $U$. If $s_0, \ldots, s_k$ give local coordinates at a point $x \in U$ and, say, $s_0(x) \neq 0$, then $s_1/s_0, \ldots, s_k/s_0$ are analytically independent meromorphic functions on $U$, hence algebraically independent. Therefore, the transcendence degree $e$ of $\mathcal{Q}(E|_U)$ equals $n$. Consequently, $\mathcal{K}(U) = \mathcal{Q}(E|_U)$ and for any given meromorphic function $\mu$ on $U$ there exists sections $s, t \in H^0(U, E^m)$, for some $m \geqslant 0$, such that $\mu = s/t$.

Assume that we have the hypothesis (H1). The manifold $M$ being $q$–complete $(q < n)$ it is known [2] that $H^{n-1}(M, F) = 0$ for any holomorphic vector bundle $F$ on $M$. By results of Hörmander [10] we have

$$H^{n-1}(M, F) \cong \mathcal{H}^{n-1}(M, F)$$

so in particular
$$\mathcal{H}^{n,n-1}(M,(E^m)^*) \cong H^{n-1}(M, \Omega^n((E^m)^*)) = 0.$$
Therefore any $CR$ section of $E^m$ over $bM$ has a holomorphic extension to $M$ by Lemma 2.1. Applying this to the traces of $s$ and $t$ on $bM$ together with the uniqueness theorem we get sections $\tilde{s}, \tilde{t} \in H^0(M, E^m)$ extending $s$ and $t$. Then $\tilde{s}/\tilde{t}$ is a meromorphic function extending $\mu$.

Assume now (H2). Theorem 2.2 shows that there exists a non-trivial holomorphic function in the neighbourhood of $M$. Therefore, for any holomorphic vector bundle $F$ Theorem 2.3 gives a line bundle $L$ such that any holomorphic section of $F \otimes L$ in the neighbourhood of $bM$ extends holomorphically in the interior of $M$. An important feature of the line bundle $L$ is that it is trivial near $bM$: from the divisor of a suitable holomorphic function on $Y$ we cut out all the components intersecting the complement of $M$. Consequently any section in $E^m$ near $bM$ may be viewed as a section in $E^m \otimes L$ and the same reasoning as above concludes the proof.

Finally we prove the conclusion assuming (H3). Since we do not have a non-constant holomorphic function anymore we have to prove a version of Theorem 2.3. For $k \geqslant 0$ let us denote $E^{-k} = (E^*)^k$.

**Lemma 2.5.** *Let $Y$ be a Moishezon manifold $Y$ of dimension $n \geqslant 2$ and $M$ a smooth domain satisfying (H3). Then for every holomorphic line bundle $F$ on $Y$ there exist $k \in \mathbb{N}$ and a holomorphic line bundle $L$ such that every smooth CR section of $F \otimes L \otimes E^{-k}$ over $bM$ extends by a smooth section of $F \otimes L \otimes E^{-k}$ over $\overline{M}$, holomorphic in $M$.*

*Proof.* Consider the map
$$H^0(Y, E^k) \times \mathcal{H}^{n,n-1}(M, F^*) \longrightarrow \mathcal{H}^{n,n-1}(M, E^k \otimes F^*),$$
(2.3)
$$(s, \sigma) \longmapsto H(s \wedge \sigma)$$
where $H$ is the projection on the space of harmonic forms and $s \wedge \sigma$ is the natural sesquilinear pairing. It suffices to show that the following claim is true:

*Claim.* There exists $k \in \mathbb{N}$ and a section $s \in H^0(Y, E^k)$ such that for any section $\alpha \in \mathcal{H}^{n,n-1}(Y, F^*)$ we have $H(s \wedge \alpha) = 0$.

Indeed, let $\sigma$ be any smooth $CR$ section of $F$ over $bM$. Then the following relation holds for every $\chi \in \mathcal{H}^{n,n-1}(M, F^*)$:
$$\int_{bM} \chi \wedge (s \otimes \sigma) = \int_{bM} (\chi \wedge s) \otimes \sigma = \int_{bM} H(\chi \wedge s) \otimes \sigma = 0.$$
Since $E$ is trivial in the neighbourhood of $bM$, $s \otimes \sigma$ is a section of $F$ over $bM$ so by Lemma 2.1 the section $s \otimes \sigma$ has a smooth extension to $\overline{M}$ which is holomorphic in $M$. Let us consider the divisor $D$ consisting of the irreducible components (counted with multiplicities) of the zero set of $s$ contained in $M$. The line bundle $L$ is obtained as the line bundle associated to $D$: $L = [D]$.

Remark that $L$ is trivial in the neighbourhood of $bM$. Let us consider a holomorphic section $\zeta$ in $L$ such that the divisor of $\zeta$ is exactly $D$, i.e., $D = (\zeta)$

Suppose that $\varsigma \in \mathcal{C}^\infty(bM, F \otimes L \otimes E^{-k})$ is such that $\bar{\partial}_b \varsigma = 0$. Since $L$ and $E$ are trivial in the neighbourhood of $bM$, $\varsigma$ induces a $CR$ section $\sigma \in \mathcal{C}^\infty(bM, F)$. Then $s \otimes \sigma \in \mathcal{C}^\infty(bM, F)$ has an extension $S \in \mathcal{C}^\infty(\overline{M}, F)$ which is holomorphic in $M$. We define a meromorphic section of $F \otimes L \otimes E^{-k}$ by
$$\Sigma = s^{-1}\zeta S.$$
Its pole set in $M$ is a union of branches of the zero set of $s$ which intersect $bM$. Since $L$ is trivial in the neighbourhood of $bM$, $\zeta = 1$ in that neighbourhood, so that the boundary value of $\Sigma$ is $\sigma = \varsigma$. We infer that the pole set of $\Sigma$ does not intersect $bM$ so it is empty. Therefore, $\Sigma$ is a holomorphic extension in $M$ of the given $\varsigma$.

We turn now to the Claim. Assume that the assertion of the claim would be false. Then a simple algebraic argument shows that

(2.4) $\qquad h^0(Y, E^k) \leqslant h^{n,n-1}(M, F^*) \times h^{n,n-1}(M, E^k \otimes F^*)$

for the dimensions of the vector spaces in (2.2). Set $n = \dim_{\mathbb{C}} Y$. Since $F$ has Kodaira dimension $\kappa(F) = n$ we have

(2.5) $\qquad h^0(Y, E^k) \approx k^n, \quad \text{for} \quad k \longrightarrow \infty$

(that is, there are constants $C, C' > 0$ such that $C'k^n \leqslant h^0(Y, E^k) \leqslant Ck^n$ for large $k$). But the right-hand member of (2.4) is a constant due to the fact that $E$ is trivial on $M$. Thus for large $k$ (2.4) cannot possibly hold. This finishes the proof of the Lemma. $\square$

Let $\mu$ be a meromorphic function on $U$ which can be represented as a quotient $\mu = s/t$ for $s, t \in H^0(U, E^m)$ for suitable $m \geqslant 0$. We apply Lemma 2.5 to $F = E^m$ so that there exist $k$ and $L$ such that the conclusion of the Lemma holds. Since $L$ and $E^{-k}$ are trivial near the boundary we can view $s$ and $t$ as sections in $F \otimes L \otimes E^{-k}$ near $bM$ and we get sections $\tilde{s}$ and $\tilde{t}$ of $F \otimes L \otimes E^{-k}$ over $M$ extending them. Then the meromorphic function $\tilde{s}/\tilde{t}$ equals $\mu$ near $bM$ and extents it over $M$. $\square$

*Remark 2.6.* A more transparent proof of the extension of sections result under hypothesis (H1) can be given in the spirit of the Ehrenpreis [6] proof of the Hartogs theorem. One may use the the cohomology groups with compact support but in order to maintain ourselves in the spirit of the work of Kohn and Rossi too let us consider the Dirichlet cohomology groups:

$$H^{0,1}_0(M, F) = \frac{\{(0,1) - \text{forms } \psi \text{ smooth on } \overline{M}, \bar{\partial}\psi = 0, \psi|_{bM} = 0\}}{\{\bar{\partial}\eta \text{ for } \eta \text{ smooth function on } \overline{M}, \eta|_{bM} = 0, \bar{\partial}\eta|_{bM} = 0\}}$$

If the Levi form of $bM$ has at least one positive eigenvalue then by the Kohn-Rossi version of the Serre duality
$$H^{0,1}_0(M, F) = H^{n,n-1}(M, F^*)$$

and under (H1) these groups vanish. If $\varphi$ is a $CR$ section of $F$ over $bM$, then there exists a smooth extension $\varphi'$ to $\overline{M}$ such that $\bar{\partial}\varphi'|_{bM} = 0$ (see the proof of 3.2). We can then solve the equation $\bar{\partial}\eta = \bar{\partial}\varphi'$ with $\eta|_{bM} = 0$ and $\varphi' - \eta$ will be the desired holomorphic extension. This argument can be carried out for forms too yielding however a weak result: for any $(0, p)$–form $\phi$ with values in $F$, $p < n - q$, one can find a $\bar{\partial}$–closed form in $M$ having the same tangential values on $bM$ as $\phi$. For the extension of $CR$ forms see the work of C. Laurent-Thiébaut and J. Leiterer [12].

*Remark 2.7.* The hypothesis (H3) can be weakened and the conclusion of the Theorem 2.4 still holds. Namely one can assume that $E \to Y$ is semi-positive with Kodaira-Iitaka maximum dimension and $E$ is trivial only in the neighbourhood of the boundary $bM$. The proof is the same; the single diffrence is that the failure of (2.4) results from asymptotic Morse inequalities [4], that is

$$(2.6) \quad \dim H^{n,n-1}(M, E^k \otimes F^*) \leqslant \frac{k^n}{n!} \int_{M(n-1)} (-1)^{n-1} \left(\frac{i}{2\pi}\mathbf{c}(E)\right)^n + o(k^n)$$

for $k \longrightarrow \infty$. In this formula $i\mathbf{c}(E)$ is the curvature of $E$ and $M(n-1)$ is the open set of $M$ where $i\mathbf{c}(E)$ has exactly 1 positive eigenvalue and $n-1$ negative ones. Since $E$ is supposed to be semi-positive on $M$ and $n \geqslant 2$ we infer that $M(n-1) = \varnothing$ and the curvature integral in (2.6) has to vanish. Therefore the right-hand side of (2.4) is $o(k^n)$ so that this inequality cannot be true for large $k$. Unfortunately the author has no example in which this situation holds, unless the bundle $E$ is trivial in a whole neighbourhood of $M$, i.e we have (H3). For the latter see Example 2.11.

Another application of Theorem 2.3 is the extension of analytic hypersurfaces across strictly pseudoconcave boundaries. It is well-known that analytic sets of dimension $\geq 2$ may be extended across such boundaries in some small neighbourhood. For analytic curves this is not always possible; see a counterexample in [9] (§1, subsection 4). That is why we shall consider the case when the ambient manifold has dimension $\geqslant 3$. Let $M$ be a strongly pseudoconvex domain, i.e., the Levi form of $bM$ is positive definite. First remark that $M$ is $1$–convex, in the sense that there exists an exhaustion function $\varrho$ of $M$ which is $1$–convex (strictly plurisubharmonic) outside a compact set $K$. For almost all $c > \sup \varrho|_K$, $M_c = \{x \in M : \varrho(x) < c\}$ is stongly pseudoconvex. The boundary of $M \setminus \overline{M}_c$ has two components: one with positive definite Levi form and one with negative definite Levi form.

**Proposition 2.8.** *Let $M$ a strongly pseudoconvex domain of dimension $\geqslant 3$. Then for any analytic hypersurface $H$ in $M \setminus \overline{M}_c$ there exists an analytic hypersurface $\widehat{H}$ in $M$ such that $H \subset \widehat{H}$.*

*Proof.* Let $E$ be the holomorphic line bundle on $M \setminus \overline{M}_c$ associated to the divisor $H$ and let $\mathcal{E}$ be the locally free sheaf of holomorphic sections in $E$.

By the reduction principle of Grauert [7], there exist a Stein complex space $\widetilde{M}$, a proper holomorphic surjection $\gamma : M \longrightarrow \widetilde{M}$, a finite subset $A \subset \widetilde{M}$ such that $\gamma^{-1}(x)$ is a connected nowhere dense analytic set of $M$ for $x \in A$ and

$$\gamma : M \setminus \gamma^{-1}(A) \longrightarrow \widetilde{M} \setminus A$$

is biholomorphic. In particular we see that there exist non-constant holomorphic functions on $M$. We may assume that $\gamma^{-1}(A) \subset K$ so that the direct image $\gamma_*(\mathcal{E})$ is a locally free analytic sheaf on $\widetilde{M} \setminus \widetilde{M}_c$, where $\widetilde{M}_c = \gamma(\overline{M}_c)$. Let us define the first absolute gap-sheaf $\gamma_*(\mathcal{E})^{[1]}$ which is the sheaf associated to the presheaf $U \longrightarrow \liminf \Gamma(U \setminus A, \gamma_*(\mathcal{E}))$ where $A$ runs over all analytic subsets of $U$ of dimension $\leqslant 1$. Because $\gamma_*(\mathcal{E})$ is locally free and the dimension of $\widetilde{M}$ is $\geqslant 3$ we see by the Riemann removable singularity theorem that $\gamma_*(\mathcal{E})^{[1]} = \gamma_*(\mathcal{E})$. Since $\widetilde{M}$ contains no positive-dimensional compact analytic set, Corollary 6.1 of Andreotti and Siu [3] shows that there exists a coherent analytic sheaf $\mathcal{F}$ on $\widetilde{Y}$ such that $\mathcal{F}|_{\widetilde{M} \setminus \widetilde{M}_c} = \gamma_*(\mathcal{E})$. Then $\gamma^*(\mathcal{F})$ is a coherent analytic sheaf on $M$ such that $\gamma^*(\mathcal{F})|_{M \setminus \overline{M}_c} \cong \mathcal{E}$. By Proposition 6.3 of Andreotti and Siu [3] we can take $\gamma^*(\mathcal{F})$ to be locally free, that is, $E$ can be extended to a holomorphic line bundle $F$ on $M$. Let us consider a section $s$ of $E$ whose divisor is $H$. By Theorem 2.3 applied to a domain $M_d$ with $d > c$ there exists a holomorphic line bundle $L$ on $M$, trivial on $bM_d$, such that every holomorphic section of $F \otimes L$ in the neighbourhood of $bM_d$ can be extended to $M_d$. Since $L$ is trivial on $bM_d$ the section $s$ can be viewed as a section in $F \otimes L$ in the neighbourhood of $bM_d$. There exists a global section section $S$ extending it and by considering the components of the divisor of $S$ we find a hypersurface $\widehat{H}$ of $M$ such that $H \subset \widehat{H}$. $\square$

*Remark 2.9.* It easily seen that the same reasoning shows that any meromorphic function defined in the neighbourhood of $bM$ of a strongly pseudoconvex domain $M$ has a meromorphic extension to $M$. We do not need neither the domain to be a subset of a Moishezon manifold nor the meromorphic function to be defined in a pseudoconcave neighbourhood of $bM$. It suffices only to remark that for any meromorphic function $\mu$ defined on a complex manifold $U$ there exist a holomorphic line bundle such that the meromorphic function can be written as the quotient of two global holomorphic sections of that bundle. To this end we choose a covering $U_i$ of $U$ such that $\mu = f_i/h_i$ on $U_i$ and the germs of the holomorphic functions $f_i, h_i$ at every point of $U_i$ are relatively prime. Since $f_i/h_i = f_j/h_j$ on $U_i \cap U_j$ we see that there exists a holomorphic function $g_{ij}$ on $U_i \cap U_j$ such that $f_i = g_{ij} f_j$ and $h_i = g_{ij} h_j$. Then $\{g_{ij}\}$ is the cocycle of the desired line bundle. This line bundle at hand we can extend it, as well as its sections, as in the previous proof and then take the quotient of the extended sections. This statement can be viewed as a particular case of the powerful result of Ivashkovich [11] already quoted in the Introduction.

Let us give a few examples of domains for which our theorem applies. We begin with an example of domain whose Levi form of the boundary has at least one positive eigenvalue and an example of pseudoconcave manifold.

**Example 2.10.** Let $Z$ be a compact complex manifold and let $Y$ be an analytic subset such that $\dim Y = q \leqslant \dim Z - 2$. Let $X = Z \setminus Y$. Then there exists a function $\varrho : X \longrightarrow \mathbb{R}$ such that

- for any $c \in \mathbb{R}$ we have $X_c = \{x \in X : \varrho(x) > c\} \Subset X$, that is, $\varrho$ is an exhaustion function from bellow;
- $\varrho$ is $(q+1)$–convex outside a compact set of $X$.

For the proof see Ohsawa [15]. Therefore, for sufficiently small $c$ the boundary of the set $X \setminus \overline{X}_c$, if smooth, has a Levi form with at least one positive eigenvalue. Consequently, for small $c$ any neighbourhood of the set $X_c$ is pseudoconcave. Note that, in the terminology of [2], $X$ is a particular case of $q$–concave manifold.

We give now an example of the situation in hypothesis (H3).

**Example 2.11.** Let $Y$ be a Moishezon manifold carrying a line bundle $E$ with maximum Kodaira-Iitaka dimension. Let $s \in H^0(Y, E)$, $D = \{s = 0\}$ so that $E$ is trivial on the complement of $D$. Let $\pi : \widetilde{Y} \longrightarrow Y$ the blow up of center $x \in Y \setminus D$. Then $\pi^{-1}(D)$ and $\pi^{-1}(x)$ are compact disjoint sets. The pull-back $\widetilde{E} = \pi^*(E)$ has maximum Kodaira-Iitaka dimension and it is trivial on the complement of $\pi^{-1}(D)$. Let us consider an analytic set $A \subset \widetilde{Y} \setminus \pi^{-1}(D)$ of dimension $q \leqslant \dim Y - 2$ (one can take $A \subset \pi^{-1}(x) \cong \mathbb{P}^{\dim Y - 1}$) so that $X = \widetilde{Y} \setminus A$ admits an exhaustion function from bellow as in the preceding example. Therefore, $\widetilde{Y}$, $\widetilde{E}$ and $M = \widetilde{Y} \setminus X_c$ for convenient very small $c$ (see 2.10) satisfy the hypothesis (H3): the Levi form of $bM$ has $n - q - 1$ positive eigenvalues and $\widetilde{E}$ is trivial on $M$ since $X_c$ contains the compact set $\pi^{-1}(D)$ for $c < \inf\{\varrho(x) : x \in \pi^{-1}(D)\}$.

We end with an example of $q$–complete manifold and of a very small pseudoconcave manifold.

**Example 2.12.** If $Z$ is a compact complex manifold and $E$ is a positive holomorphic line bundle over $Z$ and consider the global sections $s_1, \ldots, s_q$, $q \leqslant \dim Z - 1$ and $Y = \{x \in Z : s_1(x) = \cdots = s_q(x) = 0\}$. Then $X := Z \setminus Y$ is $q$–complete. Indeed, if $s_k^j$ are the local representations of these sections as holomorphic functions and if $h_j$ are the local representation of the positive metric on $E$ then $\varrho := -\log(h_j \sum_k |s_k^j|^2)$ is a $q$–convex exhaustion function for $X$ ($i\partial\bar\partial\varrho = -i\partial\bar\partial \log h_j - i\partial\bar\partial \log \sum_k |s_k^j|^2$; the first term of the sum is positive definite being the curvature of $E$ and the second vanishes on the set $\{s_2/s_1 = \cdots = s_q/s_1\}$ (say $s_1 \neq 0$) which has dimension at least $n - q + 1$ so that $i\partial\bar\partial\varrho$ has at least $n - q + 1$ positive eigenvalues). Thus $Z \setminus X_c$ is pseudoconcave for large $c$. Since $Z \setminus X_c$, $c \in \mathbb{R}$, forms a fundamental system of neighbourhoods of $Y$ every neighbourhood of $Y$ is pseudoconcave.

Using the last example we obtain the following striking result. Let $\mu$ be a meromorphic function defined in a neighbourhood of the set $Y$ above; then $\mu$ extends as a meromorphic function to $Z$ and $Z$ being projective (it admits a positive line bundle) $\mu$ is in fact rational by the well-known theorem of Weierstrass, Hurwitz and Chow. But this is only a particular case of a more general result of Dingoyan [5] asserting that a meromorphic function on a pseudoconcave projective manifold is in fact rational.

## 3. EXTENSION OF $CR$ SECTIONS FROM (WEAKLY) PSEUDOCONVEX BOUNDARIES

In this section we consider (weakly) pseudoconvex domains in complex manifolds. A smooth domain $M$ in a complex manifold $X$ is said to be pseudoconvex if the Levi form of its boundary is positive semi-definite everywhere. The extension procedure used so far involves the the global regularity of the $\bar{\partial}$–Neumann problem. That is why we consider the supplementary condition that the dual of the bundle where the section to be extended lives is positive near $bM$. This condition was introduced by Kohn [14] in the following form: there exists a smooth, non-negative function $\lambda : X \longrightarrow \mathbb{R}$ which is strictly plurisubharmonic in the neighbourhood of $bM$. This hypothesis implies that that $M$ is $1$–convex ([14], Theorem 3.7) and any line bundle is positive near $bM$. Let us introduce the following hermitian metrics on the trivial line bundle $X \times \mathbb{C}$:

$$g_t = \exp(-t\lambda), \quad t \geqslant 0.$$

This gives $L^2$ inner products $(\cdot, \cdot)_t$ on $M$, Hodge operators $\#_t$, formal adjoints $\vartheta_t = -\#_t \bar{\partial} \#_t$ and hilbertian adjoints $\bar{\partial}_t^*$ of $\bar{\partial}$. The $L^2$ spaces obtained by completing $\mathcal{E}^{p,q}(\overline{M})$ under the norms $\|\cdot\|_t$ are denoted $\mathcal{L}_t^{p,q}(M)$ (topologically all are the same). Consider the Laplace-Beltrami operator $\Delta_t'' = \bar{\partial}\bar{\partial}_t^* + \bar{\partial}_t^*\bar{\partial}$ defined by:

$$\mathrm{Dom}(\Delta_t'') = \{u \in \mathrm{Dom}(\bar{\partial}) \cap \mathrm{Dom}(\bar{\partial}_t^*) \,:\, \bar{\partial}u \in \mathrm{Dom}(\bar{\partial}_t^*),\, \bar{\partial}_t^* u \in \mathrm{Dom}(\bar{\partial})\},$$

and the kernel of $\Delta_t''$:

$$\mathcal{H}_t^{p,q}(M) = \{u \in \mathrm{Dom}(\bar{\partial}) \cap \mathrm{Dom}(\bar{\partial}_t^*) \,:\, \bar{\partial}u = 0,\, \bar{\partial}_t^* u = 0\},$$

together with the orthogonal projector $H_t$ on this subspace. For $s \in \mathbb{N}$ denote by $\mathcal{E}_s^{p,q}(\overline{M})$ the space of $(p,q)$–forms with $s$ times continuously differentiable coefficients on $\overline{M}$. The result about the global regularity is the following.

**Theorem 3.1** (Kohn [14]). *Let $M$ be a smooth pseudoconvex domain in the complex manifold $X$; assume that there exists a smooth, non-negative real function $\lambda$ on $X$ which is strictly plurisubharmonic in the neighbourhood of $bM$. Then for any $s \in \mathbb{N}$ there exists a real $T_s$ such that for $t \geqslant T_s$ and $q \geqslant 1$ we have that:*

(i) $\Delta_t''$ *has closed range and* $\mathcal{H}_t^{p,q}(M)$ *is a finite dimensional subspace of* $\mathcal{E}_s^{p,q}(\overline{M})$.

*There exists a bounded operator* $N_t : \mathcal{L}_t^{p,q}(M) \longrightarrow \mathcal{L}_t^{p,q}(M)$ *such that:*

(ii) $\mathrm{Range}(N_t) \subset \mathrm{Dom}(\Delta_t'')$ *and we have the strong Hodge decomposition:*
$$\alpha = \bar{\partial}\bar{\partial}_t^* N_t \alpha + \bar{\partial}_t^* \bar{\partial} N_t \alpha + H_t \alpha, \quad \alpha \in \mathcal{L}_t^{p,q}(M).$$

(iii) $N_t$ *commutes with* $\Delta_t''$, $\bar{\partial}$, $\bar{\partial}_t^*$ *on their domains.*

(iv) $N_t$ *and* $H_t$ *map* $\mathcal{E}^{p,q}(\overline{M})$ *into* $\mathcal{E}_s^{p,q}(\overline{M})$.

Therefore we have the following analogue of Lemma 2.1:

**Lemma 3.2.** *Let* $s \geqslant 1$ *and* $t \geqslant T_{s+1}$ *like in Theorem 3.1. Then for any smooth CR function* $\varphi$ *on* $bM$ *there exists* $\varphi_0 \in \mathcal{E}_s^{0,0}(\overline{M})$ *such that* $\varphi_0|_{bM} = \varphi$ *and* $\bar{\partial}\varphi_0 = 0$ *in* $M$ *if and only if*

$$\text{(3.1)} \qquad \int_{bM} \chi \wedge \varphi = 0, \quad \text{for all } \chi \in \mathcal{H}_t^{n,n-1}(M).$$

*Proof.* Let $\varphi'$ be an arbitrary smooth extension of $\varphi$ to $M$. Then we can achieve that $\bar{\partial}\varphi'|_{bM} = 0$ by replacing $\varphi'$ by

$$\text{(3.2)} \qquad \widehat{\varphi} = \varphi' - r\,\bar{L}\,\varphi'$$

where $r$ is a defining function for $M$ and $L$ is a $(1,0)$ vector field which satisfy $L\,r = 1$ near $bM$. As before we consider

$$\alpha = \#_t \bar{\partial}\varphi' = -\vartheta_t \#_t \varphi' \in \mathcal{E}^{n,n-1}(\overline{M}).$$

Since $\alpha|_{bM} = 0$ we have $\alpha \in \mathrm{Dom}(\bar{\partial}_t^*)$. It is clear that $\bar{\partial}_t^* \alpha = 0$ and by (iii) of Theorem 3.1 this implies that $\bar{\partial}_t^* N_t \alpha = 0$. Then (ii) shows that $\alpha = \bar{\partial}_t^* \bar{\partial} N_t \alpha$ if and only if $H_t \alpha = 0$. Since $\#_t \#_t = id$ on functions the last condition is equivalent to 3.1 by integration by parts. Consider the form $\psi = \bar{\partial} N_t \alpha$. By (iv) applied to $s+1$ we infer that $\alpha \in \mathcal{E}_s^{n,n}(\overline{M})$. Moreover $\psi \in \mathrm{Dom}(\bar{\partial}_t^*) \cap \mathcal{E}_s^{n,n}(\overline{M})$ so that $\psi|_{bM} = 0$ because $s \geqslant 1$. Define now $\varphi_0 \in \mathcal{E}_s^{0,0}(\overline{M})$ by $\#\varphi_0 = \#\varphi' + \psi$. It is now clear that $\varphi_0$ is the desired extension. □

We can state now an analogue of Theorem 2.2.

**Theorem 3.3.** *Let* $M$ *be a smooth pseudoconvex domain with connected boundary; assume that there exists a smooth, non-negative real function* $\lambda$ *on* $X$ *which is strictly plurisubharmonic in the neighbourhood of* $bM$. *Then for any smooth CR function* $\varphi$ *on* $bM$ *and any* $s \geqslant 1$ *there exists an* $s$ *times continuously differentiable extension* $\varphi_0$ *on* $\overline{M}$, *holomorphic in* $M$.

The proof is essentially the same as of Theorem 2.2 from [13]. For the convenience of the reader we sketch it here.

*Proof.* Remark first that the set of $CR$ functions forms a ring. If the only $CR$ functions on $bM$ are the constants then the conclusion follows. Assume now that there exist at least one non-constant $CR$ function $f$; it necessarily takes an infinity of values. The key observation is that we can find a monic polynomial $P$ with complex coefficients such that $P(f)$ has an extension in $\mathcal{E}_s^{0,0}(\overline{M})$, holomorphic in $M$. Indeed, take a basis $\chi_1, \ldots, \chi_k$ in $\mathcal{H}_t^{n,n-1}(M)$, for some fixed $t \geqslant T_{s+1}$, and let $c_{ij} = \int_{bM} f^i \chi_j$, $i = 0, \ldots, k$, $j = 1, \ldots, k$. Choose complex numbers $a_i$, $i = 0, \ldots, k$ such that $\sum_i a_i c_{ij} = 0$ for all $j$ and $a_k = 1$. Put $P(z) = \sum_{i=0}^k a_i z^i$. Since $\bar{\partial}_b f^i = 0$ we have $\bar{\partial}_b P(f) = 0$. Also $\int_{bM} P(f) \chi = 0$ for every $\chi \in \mathcal{H}_t^{n,n-1}(M)$. Therefore by Lemma 3.2 there exist $h \in \mathcal{E}_s^{0,0}(\overline{M})$ such that $h = P(f)$ on $bM$ and $\bar{\partial} h = 0$ in $M$. Since $P(f)$ hits an infinity of values, $h$ is non-constant.

Now the idea is to write $f$ as a quotient of two functions in $\mathcal{E}_s^{0,0}(\overline{M})$ holomorphic in $M$. Using the same reasoning as above but replacing $f^i$ with $h^i f$, $i = 0, \cdots, k$, we find a polynomial $Q$ such that $Q(h)f$ has an extension $F \in \mathcal{E}_s^{0,0}(\overline{M})$, holomorphic in $M$. Set $G = Q(h)$. Of course $f = F/G$ on $bM$ and $F/G$ is holomorphic in the complement of the zero set of $G$. If we show that $F/G$ has a holomorphic extension to all of $M$ we are done. For this we use the Riemann removable singularity theorem. Note that $h$ being non-constant so is $G = Q(h)$ to the effect that $\{G = 0\}$ is thin. So we have to show that $F/G$ is locally bounded in $M$. For this purpose we return to the polynomial $P$. First consider the holomorphic functions $G^k P(F/G)$ and $G^k h$; they have the same boundary values so they are identical. Consequently $P(F/G) = h$ on $M \setminus \{G = 0\}$. Since $h$ is locally bounded $P(F/G)$ is locally bounded. But $P$ being monic this entails that $F/G$ is locally bounded, too. □

It is known [8] that in general there exists no strictly plurisubharmonic function near the boundary $bM$. However, in Grauert's example there exists a positive line bundle near $bM$. This case was studied by K. Takegoshi [17] who extended the regularity theorem as follows.

**Theorem 3.4** ([17]). *Let $M$ be a smooth domain in a complex manifold $X$ and let $E \to X$ be a holomorphic line bundle which is positive near $bM$. Then for every $s \in \mathbb{N}$ there exists a rank $m(s)$ such that for $m \geqslant m(s)$ and $q \geqslant 1$ we have that:*

(i) $\Delta_m''$, *the Laplace-Beltrami operator acting on $\mathcal{L}^{p,q}(M, E^m)$, has closed range and its kernel $\mathcal{H}^{p,q}(M, E^m)$ is a finite dimensional subspace of $\mathcal{E}_s^{p,q}(\overline{M}, E^m)$.*

*There exists a bounded operator $N_m : \mathcal{L}^{p,q}(M, E^m) \longrightarrow \mathcal{L}^{p,q}(M, E^m)$ such that:*

(ii) $\mathrm{Range}(N_m) \subset \mathrm{Dom}(\Delta_m'')$ *and we have the strong Hodge decomposition:*
$$\alpha = \bar{\partial} \bar{\partial}_t^* N_m \alpha + \bar{\partial}_t^* \bar{\partial} N_m \alpha + H_m \alpha, \quad \alpha \in \mathcal{L}^{p,q}(M, E^m).$$

(iii) $N_m$ *commutes with $\Delta_m''$, $\bar{\partial}$, $\bar{\partial}_m^*$ on their domains.*

(iv) $N_m$ *and $H_m$, the orthogonal projection on $\mathcal{H}^{p,q}(M, E^m)$, map $\mathcal{E}^{p,q}(\overline{M}, E^m)$ into $\mathcal{E}_s^{p,q}(\overline{M}, E^m)$.*

Thus we have the following version of Lemma 3.2

**Lemma 3.5.** *Let $s \geqslant 1$ and $m \geqslant m(s+1)$ like in Theorem 3.4. Then for any smooth CR section $\varphi \in \mathcal{C}^\infty(bM, E^{-m})$ there exists $\varphi_0 \in \mathcal{E}_s^{0,0}(\overline{M}, E^{-m})$* *such that $\varphi_0|_{bM} = \varphi$ and $\bar\partial \varphi_0 = 0$ in $M$ if and only if*

$$(3.3) \qquad \int_{bM} \chi \wedge \varphi = 0, \quad \text{for all } \chi \in \mathcal{H}^{n,n-1}(M, E^m).$$

The only thing we must prove is that there exists a smooth extension $\varphi'$ such that $\bar\partial \varphi'|_{bM} = 0$. This is easily seen because formula (3.2) carries over the bundle case. Indeed, the transition functions are holomorphic and the vector field $\bar L$ of type $(0,1)$, so that $\bar L \varphi'$ is still a smooth section.

We shall now restrict attention to a particular case of pseudoconvex domains. Let $X$ be a weakly 1–complete manifold, *i.e.* there exists a plurisubharmonic smooth exhaustion function $\varrho : X \longrightarrow \mathbb{R}$. If $c$ is a regular value of $\varrho$ then $M = X_c = \{\varrho < c\}$ is a pseudoconvex domain and also a a weakly 1–complete manifold.

**Theorem 3.6.** *Let $X$ be a weakly 1–complete Kähler manifold and $c$ a regular value of the exhaustion function. Assume that $E \to X$ is a semi-positive holomorphic line bundle which is positive near $bX_c$. Then for $s \geqslant 1$ any CR section $\varphi \in \mathcal{C}^\infty(bX_c, E^{-m})$ extends to a section $\varphi_0 \in \mathcal{E}^{0,0}_s(\overline{X}_c, E^{-m})$ which is holomorphic in $X_c$.*

*Proof.* Let $m \geqslant m(s+1)$ such that Lemma 3.5 holds. By the representation theorem of [18]

$$\mathcal{H}^{n,n-1}(X_c, E^m) \cong H^{n,n-1}(X_c, E^m) \cong H^{n-1}(X_c, \Omega^n(E^m)),$$

where $n = \dim X$. But now the hypotheses we made on $X_c$ and $E$ are the same with the hypotheses of the vanishing theorem of Takegoshi [19] to the effect that these groups vanish. Therefore relation (3.3) is trivially satisfied and the conclusion follows by Lemma 3.5. □

We obtain thus an asimptotically vanishing theorem for the boundaries of some special pseudoconvex domains.

**Corollary 3.7.** *Assume that $E \to X$ is a positive line bundle on the weakly 1–complete manifold $X$ and the exhaustion function satisfies*

$$\operatorname{rank} \partial\bar\partial \exp(\varrho) < \dim X.$$

*If $c$ is a regular value of $\varrho$ then for sufficiently large $m$ any CR section $\varphi \in \mathcal{C}^\infty(bX_c, E^{-m})$ vanishes identically.*

*Proof.* We have only to remark that by Theorem 2.4 of Ohsawa [16], p. 221, a semi-negative line bundle $F$ of type $q$ (semi-negative with at least $\dim X - q + 1$ negative eigenvalues) on a weakly 1–complete Kähler manifold satisfies $H^p(X, F) = 0$ for $p \leqslant \dim X - q - r$, where $\operatorname{rank} \partial\bar\partial \exp(\varrho) \leqslant r$. In our case $F = E^{-m}$, $q = 1$ and $r = \dim X - 1$ so that $H^0(X_c, E^{-m}) = 0$ and the extension theorem above entails that $\varphi$ is the boundary value of the zero section. □

INSTITUTE OF MATHEMATICS OF THE ROMANIAN ACADEMY, PO BOX 1-764, RO-70700, BUCHAREST, ROMANIA

*Current address*: Department of Mathematics & Statistics, The University of Edinburgh, James Clerk Maxwell Building, King's Buildings, Mayfield Road, Edinburgh EH9 3JZ, UK

*E-mail address*: `george@maths.ed.ac.uk`